\documentclass[11pt]{article}
\usepackage[dvips]{graphics}
\title{\Large\bf Annular and circular rigid inclusions planted into a penny-shaped crack and factorization of triangular matrices}
\author{\bf Y.A.\ Antipov, S.M.\ Mkhitaryan\\ 
Department of Mathematics, Louisiana State University\\
Baton Rouge LA 70803, U.S.A.\\
Department of Mechanics of Elastic and Viscoelastic Bodies\\ National Academy of Sciences,
Yerevan 0019, Armenia}

\date{}

\newcommand{\supp}{\mathop{\rm supp}\nolimits}

\newcommand{\I}{\mathop{\rm Im}\nolimits}
\newcommand{\R}{\mathop{\rm Re}\nolimits}
\newcommand{\const}{\mbox{const}}

\newcommand{\bfm}[1]{\mbox{\boldmath ${#1}$}}

\newcommand{\beqa}{\begin{eqnarray}}
\newcommand{\eeqa}[1]{\label{#1}\end{eqnarray}}
\newcommand{\bequ}{\begin{equation}}
\newcommand{\eequ}[1]{\label{#1}\end{equation}}

\newcommand{\Md}{\partial}

\newcommand{\Ga}{\alpha}

\newcommand{\Gd}{\delta}

\newcommand{\Gg}{\gamma}
\newcommand{\Gc}{\chi}

\newcommand{\Gk}{\kappa}

\newcommand{\Gl}{\lambda}

\newcommand{\Gt}{\theta}

\newcommand{\Gr}{\rho}

\newcommand{\Gs}{\sigma}

\newcommand{\Go}{\omega}
\newcommand{\Gx}{\xi}
\newcommand{\Gy}{\psi}

\newcommand{\GD}{\Delta}
\newcommand{\GF}{\Phi}
\newcommand{\GG}{\Gamma}

\newcommand{\GO}{\Omega}

\newcommand{\GY}{\Psi}

\newcommand{\BGF}{\bfm\Phi}

\newcommand{\CD}{{\cal D}}

\newcommand{\CL}{{\cal L}}

\def\Bg{{\bf g}}

\def\By{{\bf y}}

\newcommand{\beq}{\begin{equation}}
\newcommand{\eeq}{\end{equation}}
\newcommand{\barr}{\begin{eqnarray}}
\newcommand{\earr}{\end{eqnarray}}
\newcommand{\beqn}{\begin{equation*}}
\newcommand{\eeqn}{\end{equation*}}
\newcommand{\barrn}{\begin{eqnarray*}}
\newcommand{\earrn}{\end{eqnarray*}}

\newcommand{\fr}{\frac}

\newcommand{\diag}{\mbox{diag}}

\begin{document}
\maketitle

\begin{abstract}

 Analytical solutions to two axisymmetric problems of a penny-shaped crack when an annulus-shaped (model 1) or a disc-shaped (model 2) rigid inclusion of arbitrary profile are embedded into the crack are derived. The problems are governed by integral equations with the  Weber--Sonin kernel on two segments. By the Mellin convolution theorem the integral equations  associated with the models 1 and 2 reduce to vector Riemann-Hilbert problems with $3\times 3$ and $2\times 2$ triangular matrix coefficients whose entries consist of meromorphic and of infinite indices exponential functions.
Canonical matrices of factorization are derived and the partial indices are computed. Exact representation formulas for the normal stress, the stress intensity factor, and the normal displacement are obtained and the results of numerical tests are reported.
 
\end{abstract}

\setcounter{equation}{0}

\section{Introduction}

Axisymmetric problems of the loading of penny-shaped cracks by normal tractions
applied at the crack faces in a homogeneous or a composite unbounded elastic body have been 
examined by many researchers including  \cite{sne1}, \cite{mos1}, \cite{mos2}, \cite{sne2}, \cite{wil}, \cite{ant1}.
Relevance to modeling of hydraulically induced fracture of resource bearing geological formations was a motivation  \cite{sel} for some of these studies.
The model problems admit an exact solution by quadratures or in a series form by a variety of methods such as Abelian operators,   the Wiener--Hopf technique,
orthogonal polynomials, and the Radon transform.  Motivated by modeling of fracture processes in composite elastic materials which are reinforced
with dilute concentrations of rigid circular inclusions Selvadurai and Singh analyzed \cite{sel} the problem of indentation of a penny-shaped crack 
by a smooth disc-shaped rigid inclusion. 
They employed Sneddon's integral representation \cite{sne1} of the general solution of the axisymmetric biharmonic equation
in terms of two arbitrary functions and then expressed the normal traction and displacement on the boundary
of the upper half-space through a single function. By the method of Abelian operators the governing triple
integral equation was reduced to a Fredholm integral equation of the second kind that was solved approximately by 
an asymptotic method. However, as it is shown in Section 2.1 of our paper, it is impossible to formulate the boundary conditions
in terms of a single function. This means that   
the governing
triple equations and the associated Fredholm equation are not equivalent to the model, and the asymptotic formula 
for the stress intensity factor is incorrect.

The goal of this paper is to derive an analytical solution to two problems of a penny shaped crack when an annular
(model 1) or a circular inclusion (model 2) is embedded into the crack. The inclusions are assumed to be rigid and not necessary flat.
We do not employ the theory of Abelian operators and do not end up with Fredholm integral equations. Instead, we set the problem
as an integral equation with the Weber--Sonin kernel on two segments, apply the Mellin convolution theorem and deduce
an order-3 (model 1) or order-2 (model 2) vector Riemann-Hilbert problem with a triangular matrix coefficient. To solve these problems, we 
advance the technique proposed by one of the authors \cite{ant2} for a contact model of an annulus-shaped punch. The method bypasses matrix factorization and 
eventually delivers an analytical solution that contains some
series whose coefficients solve an infinite system of linear algebraic equations. What is remarkable is that the rate of convergence 
of an approximate solution to the exact one is exponential, and when the inclusion is flat, the solution is free of quadratures.
In the case of a disc-shaped inclusion planted into a penny-shaped crack we show how the infinite system can be solved exactly in terms of recurrence relations. The same procedure is applicable
in the case of model 1 as well.
In addition to this approach by advancing further the method \cite {ant3} we factorize the $3\times 3$ and $2\times 2$ triangular matrices associated with the models. 
We prove that the factorization matrices found are not singular in any finite part of the complex plane, have the normal form
at the infinite point and constitute the canonical factorization. By analyzing the canonical factorization matrices at infinity 
we show that all the partial indices of factorization equal zero for both models.
Finally, for model 2, we derive representation formulas for the normal stress, the stress intensity factor, and the normal displacement.
Based on the exact formula for the stress intensity factor and the recurrence relations we  also obtain a simple asymptotic formula for the 
stress intensity factor.
The model we aim to analyze is an axisymmetric analog for a homogeneous space of the  two-dimensional problem  \cite{ant} concerning a rigid inclusion embedded into an interfacial crack.

\setcounter{equation}{0}

\section{Interaction of an annular inclusion and a penny-shaped crack}\label{s2}

In this section we model contact interaction of a penny-shaped crack and an annular rigid inclusion, reduce it first
to a convolution integral equation in two segments and then to a vector Riemann-Hilbert problem with a triangular $3\times 3$ matrix
coefficient. 

\subsection{Formulation}\label{s2.1}

The problem under consideration is axisymmetric one of contact interaction of a penny-shaped crack
$\{0\le r\le a, 0\le\Gt\le2\pi\}$ in the plane $z=0$ and an annular rigid inclusion $\{c\le r\le b, 0\le\Gt\le2\pi, z=\pm w(r)\}$ planted between 
the upper and lower crack faces, $0<c<b<a$. The surrounding matrix is an infinite elastic solid
whose shear modulus is $G$ and the Poisson ratio is $\nu$.
The function $w(r)$ is positive, convex, continuously differentiable, and $w(r)<<a$ everywhere in the interval $c\le r\le b$.
The inclusion surfaces are assumed to be smooth such that the tangential traction component vanishes everywhere in the contact zone 
$c_1\le r\le b_1$ ($c\le c_1<b_1\le b$). In general, the contact zone parameters $c_1$ and $b_1$ are unknown {\it a priori} and to be determined from the conditions
of boundedness of the normal contact stress
$\Gs_z$ at the points $r=c_1$ and $r=b_1$. In the particular case, when $w(r)=\Gd$, $c\le r\le b$, the inclusion is in full contact with the crack surfaces,
and $c_1=c$, $b_1=b$. 

Due to the symmetry of the problem with respect to the plane $z=0$, after the boundary conditions are linearized, it suffices to analyze the problem of the upper half-space 
with the boundary conditions in the plane $z=0$ taking the form
$$
u_z(r,0)=\left\{
\begin{array}{cc}
w(r), & c_1<r<b_1,\\
0, & r>a,\\
\end{array}
\right.
$$
\beq
\tau_{rz}(r,0)=0, \quad 0<r<\infty,
\quad \Gs_z(r,0)=0, \quad r\in(0,c_1)\cup (b_1,a).
\label{2.1}
\eeq
The elastic displacements  
and stresses may be expressed through the Love stress potential  $\Psi(r,z)$ of  the axisymmetric model by
$$
2Gu_r=-\fr{\Md^2\Psi}{\Md r\Md z}, \quad 2Gu_z=\left[2(1-\nu)\GD-\fr{\Md^2}{\Md z^2}\right]\Psi,
$$
\beq
\tau_{rz}=\fr{\Md}{\Md r}\left[(1-\nu)\GD-\fr{\Md^2}{\Md z^2}\right]\Psi,\quad
\Gs_{z}=\fr{\Md}{\Md z}\left[(2-\nu)\GD-\fr{\Md^2}{\Md z^2}\right]\Psi,
\label{2.2}
\eeq
where
\beq
\GD=\fr{\Md^2}{\Md r^2}+\fr{1}{r}\fr{\Md}{\Md r}+\fr{\Md^2}{\Md z^2}.
\label{2.3}
\eeq
The model under consideration is thus governed by the boundary value problem (\ref{2.1}) to (\ref{2.3}) for the biharmonic axisymmetric operator
\beq
\GD^2 \Psi(r,z)=0, \quad 0<r<\infty,\quad 0<z<\infty.
\label{2.4}
\eeq
As $z\to\infty$ and $0\le r<\infty$, the function $\Psi(r,z)$ and all its derivatives up to the fourth order vanish. 
The general solution to this equation is given by \cite{sne1} 
\beq
\GY(r,z)=\int_0^\infty[A(\Gx)+zB(\Gx)]e^{-\Gx z} J_0(r\Gx)d\Gx.
\label{2.4'}
\eeq
By using (\ref{2.4'}) and (\ref{2.2}) it is  verified that 
$$
\Gs_z(r,0)=\int_0^\infty[\Gx A(\Gx)+(1-2\nu)B(\Gx)]J_0(r\Gx)\Gx^2d\Gx,
$$
\beq
2Gu_z(r,0)=-\int_0^\infty[\Gx A(\Gx)+2(1-2\nu)B(\Gx)]J_0(r\Gx)\Gx d\Gx.
\label{2.4''}
\eeq
It becomes evident that there is no single function, $R(\Gx)$, which may serve in the integral representations
(\ref{2.4''}) instead of the two functions $A(\Gx)$ and $B(\Gx)$. Therefore, the triple integral equations (10) to (12) in \cite{sel}
are incorrect.

\subsection{Derivation of an order-3 vector Riemann--Hilbert problem}\label{s2.2}

To pursue our goal to reformulate the boundary value problem 
 (\ref{2.4}) as a vector Riemann--Hilbert problem, we introduce new unknown functions, $\Gc_0$, $\Gc_1$, $\psi_0$, and $\psi_1$,
and write down the first and third boundary conditions (\ref{2.1})  in the whole plane $z=0$ as
\beq
u_z(r,0)=\left\{
\begin{array}{cc}
-a\Gt_1\Gc_0(r), & 0\le r\le c_1,\\
w(r), & c_1<r<b_1,\\
-a\Gt_1\Gc_1(r), & b_1\le r\le a,\\
0, & r>a,\\
\end{array}
\right.
\quad
\Gs_z(r,0)=\left\{
\begin{array}{cc}
0, & 0\le r\le c_1,\\
\psi_0(r), & c_1<r<b_1,\\
0, & b_1\le r\le a,\\
\psi_1(r), & r>a.\\
\end{array}
\right.
\label{2.5}
\eeq
where $\Gt_1=(1-\nu)/G$. On applying the Hankel transform 
\beq
\Psi_\Gl(z)=\int_0^\infty \Psi(r,z)J_0(\Gl r)rdr,
\quad
\Gs_{z\Gl}=\int_0^\infty \Gs_z(r,0)J_0(\Gl r)rdr
\label{2.6}
\eeq
to the boundary value problem (\ref{2.4}) we deduce
$$
\left(\fr{d^4}{dz^4}-2\Gl^2\fr{d^2}{dz^2}+\Gl^4\right)\Psi_\Gl(z)=0, \quad 0<z<\infty,
  $$
 \beq
 \Gl^2(1-\nu)\Psi_\Gl(0)+\nu\fr{d^2}{dz^2}\Psi_\Gl(0) =0,
 \quad
 -\Gl^2(2-\nu)\fr{d}{dz}\Psi_\Gl(0)+(1-\nu) \fr{d^3}{dz^3}\Psi_\Gl(0)=\Gs_{z\Gl}.
 \label{2.7}
 \eeq
After the general solution to this one-dimensional boundary value problem is written down we invert
the Hankel transform and express the displacement $u_z$ in the plane $z=0$ through the normal traction as
\beq
u_z(r,0)=-\Gt_1\int_0^\infty W_{00}(r,\Gr)\Gs_z(\Gr,0)\Gr d\Gr,
\label{2.8}
\eeq
where
 $W_{mn}(r,\Gr)$ is the Weber--Sonin integral
\beq
W_{mn}(r,\Gr)=\int_0^\infty J_m(r\Gx)J_n(\Gr\Gx)d\Gx.
\label{2.10}
\eeq
Returning now to the first boundary condition in (\ref{2.1}) and using (\ref{2.5}) and (\ref{2.8}) we reformulate it as an integral equation in two segments  
\beq
\int_{c_1}^{b_1}W_{00}(r,\Gr)\psi_0(\Gr)\Gr d\Gr+\int_{a}^{\infty}W_{00}(r,\Gr)\psi_1(\Gr)\Gr d\Gr=\left\{
\begin{array}{cc}
-\Gt_1^{-1}w(r), & c_1<r<b_1,\\
0, & r>a.\\
\end{array}
\right.
\label{2.11}
\eeq
The integral equation can be recast by employing the functions  $\Gc_0(r)$ and $\Gc_1(r)$ introduced in (\ref{2.5}) and
a function $w_0(r)$ that is $w_0(r)=w(r)$, $c_1<r<b_1$
and $w_0(r)=0$ otherwise.
Extend the definitions of the functions $\Gc_j(r)$, and $\psi_j(r)$ to the whole ray $r\ge 0$ by
\beq
 \supp \Gc_0(r)\subset[0,c_1], \; \supp \Gc_1(r)\subset[b_1,a], 
\;
\supp \Gy_0(r)\subset[c_1,b_1], \; \supp \Gy_1(r)\subset[a,\infty]. 
\label{2.12}
\eeq
This brings us to the following  Mellin convolution integral equation:
\beq
\int_0^\infty l\left(\fr{r}{\Gr}\right)[\psi_0(a\Gr)+\psi_1(a\Gr)]d\Gr=
\Gc_0(ar)+\Gc_1(ar)-\fr{w_0(ar)}{a\Gt_1}, \quad 0<r<\infty,
\label{2.13}
\eeq
where
\beq
l(t)=\int_0^\infty J_0(t\Gx)J_0(\Gx)d\Gx.
\label{2.14}
\eeq
Our next step is to introduce the Mellin transforms of the functions $w(r)$, $\Gc_j(r)$, and $\psi_j(r)$ which, on account of (\ref{2.12}), are
$$
\GF_1^-(s)=\int_{\Gl_1}^1\Gc_1(ar)r^{s-1}dr,\quad
\GF_1^+(s)=\int_{1}^{1/\Gl_1}\Gc_1(b_1r)r^{s-1}dr,
$$$$
\GF_2^-(s)=\int_{\Gl_0/\Gl_1}^1\Gy_0(b_1r)r^{s}dr,
\quad
\GF_2^+(s)=\int_{1}^{\Gl_1/\Gl_0}\Gy_0(c_1r)r^{s}dr,
$$
\beq
\GF_3^-(s)=\int_{0}^1\Gc_0(c_1r)r^{s-1}dr,\quad 
\GF_3^+(s)=\int_{1}^{\infty}\Gy_1(ar)r^{s}dr,
\label{2.15}
\eeq
and evaluate the Mellin transform of the kernel $l(t)$
\beq
L(s)=\int_0^\infty l(t)t^{s-1}dt.
\label{2.16}
\eeq
Here,
\beq
 \Gl_0=\fr{c_1}{a}, \quad  \Gl_1=\fr{b_1}{a}, \quad 0<\Gl_0<\Gl_1<1.
 \label{2.17}
 \eeq 
By making use of the table integral  6.561(14) \cite{gra}
\beq
\int_0^\infty x^\mu J_\nu(ax)dx=
\fr{2^\mu \GG(1/2+\nu/2+\mu/2)}{a^{\mu+1}\GG(1/2+\nu/2-\mu/2)}, \quad
-\R \nu-1<\R\mu<\fr12, \quad a>0,
\label{2.18}
\eeq
we have
\beq
L(s)=\fr{\GG(s/2)\GG(1/2-s/2)}{2\GG(1-s/2)\GG(1/2+s/2)}, \quad 0<\R s<1.
\label{2.19}
\eeq
The functions $\Gc_0(c_1r)$ and $\psi_1(ar)$ are sought in the class of functions
having the asymptotics
\beq
\Gc_0(c_1r)=O(1), \quad r\to 0, \quad \psi_1(ar)=O(r^{-1-\Ga}), \quad r\to\infty, \quad 0<\Ga\le 1.
\label{2.20}
\eeq
Due to the Abelian theorems for the Mellin transform we conclude that the functions
$\GF_3^-(s)$ and $\GF_3^+(s)$ are analytic in the half-planes $\R s>0$
 and $\R s<\Ga$, respectively. Notice that the other functions, $\GF_1^\pm(s)$, $\GF_2^\pm(s)$, and the Mellin transform of the function 
 $w_0(ar)$, are entire functions, and therefore the Mellin transforms of all the functions under consideration are analytic at least in the strip
 $0<\R s<\Ga$.
 
 Apply now the Mellin transform to equation (\ref{2.13}). In view of the Mellin convolution theorem,
 we have the following vector Riemann--Hilbert problem with a triangular matrix coefficient:
 \beq
 \BGF^+(s)=G(s)\BGF^-(s)+\Bg(s), \quad s\in\CL,
 \label{2.21}
 \eeq
 where $\CL=\{\R s=\Gg, -\infty<\I s<+\infty\}$, $0<\Gg<\Ga\le 1$,
 $$
 G(s)=\left(
 \begin{array}{ccc}
 \Gl_1^{-s} & 0 & 0\\
0   & (\Gl_0/\Gl_1)^{-s-1} & 0\\
1/L(s)  & -\Gl_1^{s+1}  &  \Gl_0^s/L(s) \\
\end{array}
\right),
$$
\beq
\Bg(s)=\left(
 \begin{array}{c}
0\\
0\\
-\Gl_1^s[a\Gt_1 L(s)]^{-1}\hat w^-(s) \\
\end{array}
\right),\quad \hat w^-(s)=\int_{\Gl_0/\Gl_1}^1 w(b_1r)r^{s-1}dr.
\label{2.22}
\eeq
The column-vectors $\BGF^\pm(s)=(\GF_1^\pm(s),\GF_2^\pm(s),\GF_3^\pm(s))^T$ are analytic in the half-planes $\CD^\pm$,
and $\CD^+=\{\R s\le\Gg\}$,  $\CD^-=\{\R s\ge\Gg\}$.

\subsection{Solution of the vector Riemann--Hilbert problem}\label{s2.3}

Before proceeding with the solution, we note that although the matrix coefficient is a lower triangular matrix,
it is not reducible to a sequently solvable scalar Riemann--Hilbert problems. This is because the first
two problems have plus-infinite indices, and an infinite number of solutions expressible through free entire functions of
certain properties exist, 
while 
the index of the third problem is equal to $-\infty$; its solvability condition gives rise to integral equations with respect to the
entire functions coming from the first two problems \cite{ant4}. These integral equations are not simpler than
the original vector Riemann--Hilbert problem.

To derive an efficient solution to the problem (\ref{2.21}), we advance the method introduced  in  \cite{ant2}. First, we  factorize the function $L(s)$,
$$
L(s)=\fr{L^+(s)}{2L^-(s)}, \quad s\in\CL,
$$
\beq
L^+(s)=\fr{\GG(1/2-s/2)}{\GG(1-s/2)}, \quad  L^-(s)=\fr{\GG(1/2+s/2)}{\GG(s/2)}, 
\label{3.1}
\eeq
and then rewrite the third equation in (\ref{2.21}) as
\beq
\fr12L^+(s)\GF_3^+(s)=L^-(s)\GF_1^-(s)-\fr12\Gl_1^{s+1}L^+(s)\GF_2^-(s)+\Gl_0^sL^-(s)\GF_3^-(s)-\fr{\Gl_1^s}{a\Gt_1}L^-(s)\hat w^-(s).
\label{3.2}
\eeq
Next, we multiply the third equation in (\ref{2.21})  by $\Gl_1^{-s}$ and use the first equation in (\ref{2.21}) that is $\Gl_1^{-s}\GF_1^-(s)=\GF_1^+(s)$.
After rearrangement, we have
\beq
\fr{2\GF_1^+(s)}{L^+(s)}-\fr{\Gl_1^{-s}\GF_3^+(s)}{L^-(s)}=\fr{\Gl_1\GF_2^-(s)}{L^-(s)}-\left(\fr{\Gl_0}{\Gl_1}\right)^s\fr{2\GF_3^-(s)}{L^+(s)}+\fr{2\hat w^-(s)}{a\Gt_1 L^+(s)}.
\label{3.3}
\eeq
The third equation of the new system is obtained by multiplying  the third equation in (\ref{2.21})  by $\Gl_0^{-s}$. In view of the second equation in  
 (\ref{2.21}), we have 
\beq
\fr{\Gl_0^{-s}}{2}L^+(s)\GF_3^+(s)+\fr{\Gl_0}{2}L^+(s)\GF_2^+(s)-\left(\fr{\Gl_0}{\Gl_1}\right)^{-s}L^-(s)\GF_1^+(s)
=L^-(s)\GF_3^-(s)-\fr{L^-(s)}{a\Gt_1}\hat w^+(s),
\label{3.4}
\eeq
where
\beq
\hat w^+(s)=\int_1^{\Gl_1/\Gl_0}w(c_1r)r^{s-1}dr.
\label{3.5}
\eeq

Now, in the half-plane $\CD^-$, the functions $L^+(s)$ and $1/L^+(s)$ have simple poles at the points $s=2n+1$ and $s=2n+2$ ($n=0,1,\ldots$), respectively.
In the domain $\CD^+$, the function $L^-(s)$ has simple poles  at the points $s=-2n-1$, while the simple poles of $1/L^-(s)$
are $s=-2n$ ($n=0,1,\ldots$).
To remove these poles in equations (\ref{3.2}) to (\ref{3.4}), we introduce the following functions:
$$
\Psi^+(s)=\sum_{m=0}^\infty\fr{A_m^+}{s-2m-1},\quad \GO^+(s)=\sum_{m=0}^\infty\fr{B_m^+}{s-2m-2},
$$
\beq
\Psi^-(s)=\sum_{m=0}^\infty\fr{A_m^-}{s+2m+1},\quad \GO^-(s)=\sum_{m=0}^\infty\fr{B_m^-}{s+2m},
\label{3.6}
\eeq
with the coefficients $A_m^\pm$ and $B_m^\pm$ to be determined. 

We shall also need the representations
\beq
\fr{2\hat w^-(s)}{a\Gt_1 L^+(s)}=\Go_1^+(s)-\Go^-_1(s), \quad 
\fr{\hat w^+(s)}{a\Gt_1} L^-(s)=\Go_2^+(s)-\Go^-_2(s),\quad s\in\CL.
\label{3.7}
\eeq
Here, $\Go_j^\pm(s)$ are the limit values of the Cauchy inegrals
\beq
\Go_1(s)=\fr{1}{\pi i a\Gt_1}\int_{\CL}\fr{\hat w^-(\tau)d\tau}{L^+(\tau)(\tau-s)},
\quad
\Go_2(s)=\fr{1}{2\pi i a\Gt_1}\int_{\CL}\fr{\hat w^+(\tau)L^-(\tau)d\tau}{\tau-s}, 
\label{3.8}
\eeq
in the left- and right-hand sides of the contour $\CL$, respectively.

On subtracting from the left- and right-hand sides of equations (\ref{3.2}), (\ref{3.3}), and (\ref{3.4})
the functions $\Psi^+(s)$, $\GO^+(s)+\GO^-(s)$, and   $\Psi^-(s)$, respectively, using the relations (\ref{3.7}), the continuity
principle, the Liuoville theorem, and the asymptotics
$$
L^\pm(s)\sim \left(\mp\fr{s}{2}\right)^{\mp 1/2}, \quad s\in\CD^\pm, \quad s\to\infty,
$$$$
\GF_1^\pm(s)=O(s^{-3/2}), \quad \GF_2^\pm(s)=O(s^{-1/2}), \quad  s\in\CD^\pm, \quad s\to\infty,
$$
\beq
\GF_3^+(s)=O(s^{-1/2}), \quad s\in\CD^+,\quad \GF_3^-(s)=O(s^{-3/2}), \quad s\in\CD^-,\quad s\to\infty,
\label{3.9}
\eeq
we deduce the following formulas for the solution to the vector Riemann--Hilbert problem (\ref{2.21}):
$$
\GF_1^+(s)=\fr12L^+(s)[\GO^+(s)+\GO^-(s)+\Go_1^+(s)]+\fr{\Gl_1^{-s} \Psi^+(s)}{L^-(s)},
$$ 
$$
\GF_1^-(s)=\fr{\Gl_1^s}{2}L^+(s)[\GO^+(s)+\GO^-(s)+\Go_1^-(s)]+\fr{\Psi^+(s)}{L^-(s)}+\fr{\Gl_1^s \hat w^-(s)}{a\Gt_1},
$$ 
$$
\GF_2^+(s)=\fr{2}{\Gl_0}\fr{\GY^-(s)-\Go_2^+(s)}{L^+(s)}+\left(\fr{\Gl_0}{\Gl_1}\right)^{-s}\fr{L^-(s)}{\Gl_0}[\GO^-(s)+\GO^+(s)+\Go_1^+(s)],
$$
$$
\GF_2^-(s)=\fr{L^-(s)}{\Gl_1}[\GO^-(s)+\GO^+(s)+\Go_1^-(s)]+2\left(\fr{\Gl_0}{\Gl_1}\right)^{s}\fr{\Psi^-(s)-\Go_2^-(s)}{\Gl_1L^+(s)},
$$
\beq
\GF_3^+(s)=\fr{2\Psi^+(s)}{L^+(s)}, \quad \GF_3^-(s)=\fr{\Psi^-(s)-\Go_2^-(s)}{L^-(s)}.
\label{3.10}
\eeq
It is immediately seen that $\GF_1^+(s)=\Gl_1^{-s}\GF_1^-(s)$ and $\GF_2^+(s)=(\Gl_0/\Gl_1)^{-s-1}\GF_2^-(s)$.  

In general, for arbitrary selected coefficients $A^\pm_m$ and $B_m^\pm$, the functions
$\GF_1^\pm(s)$ and $\GF_2^\pm(s)$ have inadmissible simple poles. They become removable singularities if and only if the following conditions
are satisfied:
 $$
  \mathop{\rm res}\limits_{s=-2n}\GF_1^+(s)=0,\quad   \mathop{\rm res}\limits_{s=2n+1}\GF_1^-(s)=0,\quad n=0,1,\ldots,
$$
\beq
  \mathop{\rm res}\limits_{s=-2n-1}\GF_2^+(s)=0,\quad   \mathop{\rm res}\limits_{s=2n+2}\GF_2^-(s)=0,\quad n=0,1,\ldots,
\label{3.11}
\eeq 
Note that the functions $L^-(s)\GO^-(s)$ and $L^+(s)\GO^+(s)$ have removable singularities at the points $s=-2n$ and $s=2n+2$, respectively ($n=0,1,\ldots$).
The conditions (\ref{3.11}) give rise to the infinite system of linear algebraic equations with respect to $A_n^\pm$ and $B_n^\pm$
$$
B_n^-=\fr{2\Gl_1^{2n}}{\pi}\sum_{m=0}^\infty\fr{A_m^+}{n+m+1/2},
$$$$
A_n^+=\fr{\Gl_1^{2n+1}}{2\pi}
\left[\sum_{m=0}^\infty\fr{B_m^-}{n+m+1/2}+\sum_{m=0}^\infty\fr{B_m^+}{n-m-1/2}+2\Go_1^-(2n+1)\right],
$$ 
$$
A_n^-=-\fr{1}{2\pi}\left(\fr{\Gl_0}{\Gl_1}\right)^{2n+1}
\left[\sum_{m=0}^\infty\fr{B_m^+}{n+m+3/2}+\sum_{m=0}^\infty\fr{B_m^-}{n-m+1/2}-2\Go_1^+(-2n-1)\right],
$$ 
\beq
B_n^+=-\fr{2}{\pi}\left(\fr{\Gl_0}{\Gl_1}\right)^{2n+2}\left[
\sum_{m=0}^\infty\fr{A_m^-}{n+m+3/2}-2\Go_2^-(2n+2)\right], \quad n=0,1,\ldots.
\label{3.12}
\eeq
This system can be solved by the method of reduction (the rate of convergence of an approximate solution to the exact one is exponential).
Because of its structure, the system may also be solved in terms of recurrence relations. This procedure will be described in the case
of a circular inclusion in the next section.

To conclude this section, we simplify the formulas for the functions $\Go_1^\pm(s)$ and $\Go_2^\pm(s)$ 
in the case when the annular inclusion is flat. In this case $b_1=b, c_1=c$, 
$w(r)=\Gd=\const$,  and the function $\hat w^-(s)$ is simplified to the form
\beq
\hat w^-(s)=\fr{\Gd}{s}\left[1-\left(\fr{\Gl_0}{\Gl_1}\right)^s\right].
\label{3.13}
\eeq
The  integral (\ref{3.8}) can be evaluated explicitly, and the functions
$\Go_1^\pm(s)$ and $\Go_2^\pm(s)$ are written in the form
$$
\Go_1^+(s)=\fr{\Gd}{a\Gt_1}\left[\fr{2}{s}\left(\fr{1}{L^+(s)}-\fr{1}{\sqrt{\pi}}\right)-\tilde\Go^+(s)\right],
$$$$
\Go_1^-(s)=\fr{\Gd}{a\Gt_1}\left[-\fr{2}{s\sqrt{\pi}}+\fr{2}{sL^+(s)}\left(\fr{\Gl_0}{\Gl_1}\right)^s-\tilde\Go^+(s)\right],
$$
\beq
\Go_2^+(s)=\fr{\Gd}{a\Gt_1}\left[\left(\fr{\Gl_0}{\Gl_1}\right)^{-s}
\fr{L^-(s)}{s}-\tilde\Go^-(s)\right],\quad \Go_2^-(s)=\fr{\Gd}{a\Gt_1}\left[\fr{L^-(s)}{s}-\tilde\Go^-(s)\right].
\label{3.14}
\eeq
Here,
$$
\tilde\Go^+(s)=\fr{2}{\pi}\sum_{n=0}^\infty
\fr{\GG(n+3/2)}{(n+1)!}\left(\fr{\Gl_0}{\Gl_1}\right)^{2n+2}\fr{1}{s-2n-2},
$$
\beq
\tilde\Go^-(s)=\fr{1}{\pi}\sum_{n=0}^\infty
\fr{\GG(n+1/2)}{n!}\left(\fr{\Gl_0}{\Gl_1}\right)^{2n+1}\fr{1}{s+2n+1}.
\label{3.15}
\eeq
Equivalently, in terms of the hypergeometric function, these functions may be represented as
$$
\tilde\Go^+(s)=-\fr{2}{\sqrt{\pi} s}+\fr{2}{\sqrt{\pi} s}F\left(-\fr{s}{2},\fr12; 1-\fr{s}{2}; \fr{\Gl_0^2}{\Gl_1^2}\right),
$$
\beq
\tilde\Go^-(s)=\fr{\Gl_0}{\sqrt{\pi} \Gl_1(s+1)}F\left(\fr{s+1}{2},\fr12; \fr{s+3}{2}; \fr{\Gl_0^2}{\Gl_1^2}\right).
\label{3.16}
\eeq

\setcounter{equation}{0}

\section{A circular inclusion embedded into a penny-shaped crack}\label{s4}

In this section we shall examine the particular case $c=0$ of the previous model that is 
the contact interaction of a circular inclusion $\{0\le r\le b, 0\le \Gt\le 2\pi, z=\pm w(r)\}$
and a penny-shaped crack   $\{0\le r\le a, 0\le \Gt\le 2\pi\}$  in the plane $z=0$ when $b<a$.

In the notations of Section 2, we may write the governing integral equation of the problem as 
\beq
\int_{0}^{b_1}W_{00}(r,\Gr)\psi_0(\Gr)\Gr d\Gr+\int_{a}^{\infty}W_{00}(r,\Gr)\psi_1(\Gr)\Gr d\Gr=\left\{
\begin{array}{cc}
-\Gt_1^{-1}w(r), & 0<r<b_1,\\
0, & r>a.\\
\end{array}
\right.
\label{4.1}
\eeq
As before, we write the integral equation in the Mellin convolution form
\beq
\int_0^\infty l\left(\fr{r}{\Gr}\right)[\psi_0(a\Gr)+\psi_1(a\Gr)]d\Gr=
\Gc_1(ar)-\fr{w_0(ar)}{a\Gt_1}, \quad 0<r<\infty,
\label{4.2}
\eeq
where $w_0(r)=w(r)$ if $0\le r\le b_1$ and 0 otherwise, $\Gc_1(r)=0$ if $r\in[0,b_1)\cup[a,\infty)$, $\psi_0(r)=0$ if
$r>b_1$, and $\psi_1(r)=0$ if $0\le r<a$.

In the case under consideration, $\Gl_0=0$, $\Gl_1=\Gl=b_1/a\in(0,1)$, and the analogs of the Mellin transforms (\ref{2.15}) become  
$$
\GF_1^-(s)=\int_{\Gl}^1\Gc_1(ar)r^{s-1}dr,\quad 
\GF_1^+(s)=\int_{1}^{1/\Gl}\Gc_1(b_1r)r^{s-1}dr,
$$ 
\beq
\GF_2^-(s)=\int_{0}^1\Gy_0(b_1r)r^{s}dr,\quad 
\GF_2^+(s)=\int_{1}^{\infty}\Gy_1(ar)r^{s}dr.
\label{4.3}
\eeq
Due to the absence of the function $\Gc_0(r)$ and its Mellin transform, the Riemann--Hilbert problem
is now of order-2 and has the form
$$
\GF_1^+(s)=\Gl^{-s}\GF_1^-(s),
$$
\beq
\GF_2^+(s)=\fr{\GF_1^-(s)}{L(s)}-\Gl^{s+1}\GF_2^-(s)-\fr{\Gl^s\hat w^-(s)}{a\Gt_1 L(s)}, \quad s\in\CL.
\label{4.4}
\eeq
Similarly to the previous section, it can be transformed to the system of two equations
$$
\fr12L^+(s)\GF_2^+(s)=L^-(s)\GF_1^-(s)-\fr{\Gl^{s+1}}{2}L^+(s)\GF_2^-(s)-\fr{\Gl^s}{a\Gt_1}L^-(s)\hat w^-(s),
$$
\beq
\fr{\Gl\GF_2^-(s)}{L^-(s)}+\fr{2\hat w^-(s)}{a\Gt_1 L^+(s)}=\fr{2\GF_1^+(s)}{L^+(s)}-\fr{\Gl^{-s}\GF_2^+(s)}{L^-(s)},\quad s\in\CL.
\label{4.5}
\eeq
Our next step is to remove the inadmissible poles of the functions $L^+(s)$ and $1/L^-(s)$ in the right-hand sides of equations (\ref{4.5}),
use the functions $\Psi^+(s)$ and $\GO^-(s)$ introduced in 
(\ref{3.6}), the first relation in (\ref{3.7}), the continuity principle, and the Liouville theorem. This yields
$$
\fr{L^+(s)}{2}\GF_2^+(s)-\Psi^+(s)=L^-(s)\GF_1^-(s)-\fr{\Gl^{s+1}}{2}L^+(s)\GF_2^-(s)-\fr{\Gl^s}{a\Gt_1}L^-(s)\hat w^-(s)-\Psi^+(s)=0,
$$
$$
\fr{\Gl\GF_2^-(s)}{L^-(s)}+\fr{2\hat w^-(s)}{a\Gt_1 L^+(s)}-\GO^-(s)-\Go_1^-(s)=\fr{2\GF_1^+(s)}{L^+(s)}-\fr{\Gl^{-s}\GF_2^+(s)}{L^-(s)}-\GO^-(s)-\Go_1^+(s)=0,
$$
\beq
s\in D^+\cup\CL\cup\CD^-.
\label{4.5'}
\eeq
From here, we derive the solution to the vector Riemann--Hilbert problem
$$
\GF_1^+(s)=\fr12L^+(s)[\GO^-(s)+\Go_1^+(s)]+\fr{\Gl^{-s} \Psi^+(s)}{L^-(s)},
$$ 
$$
\GF_1^-(s)=\fr{\Gl^s}{2}L^+(s)[\GO^-(s)+\Go_1^-(s)]+\fr{\Psi^+(s)}{L^-(s)}+\fr{\Gl^s \hat w^-(s)}{a\Gt_1},
$$ 
\beq
\GF_2^+(s)=\fr{2\Psi^+(s)}{L^+(s)}, \quad \GF_2^-(s)=\fr{L^-(s)}{\Gl}[\GO^-(s)+\Go_1^-(s)].
\label{4.6}
\eeq
The conditions which transform the undesired simple poles of the functions $\GF^+_1(s)$ and $\GF_1^-(s)$
into removable singular points become
$$
B_n^-=\fr{2\Gl^{2n}}{\pi}\sum_{m=0}^\infty\fr{A_m^+}{n+m+1/2},
$$
\beq
A_n^+=\fr{\Gl^{2n+1}}{2\pi}
\left[\sum_{m=0}^\infty\fr{B_m^-}{n+m+1/2}+2\Go_1^-(2n+1)\right], \quad n=0,1,\ldots.
\label{4.7}
\eeq
These equations constitute an infinite system of linear algebraic equations. As in the previous section, it can be solved
numerically by the reduction method. Alternatively, its solution may be derived  in terms of recurrence relations.
For simplicity, we suppose that the inclusion is flat, $w(r)=\Gd=\const$, $0\le r\le b$. Then we have $b_1=b$ and
\beq
\hat w^-(s)=\fr{\Gd}{s}, \quad \Go_1^+(s)=\fr{2\Gd}{a\Gt_1 s}\left(\fr{1}{L^+(s)}-\fr{1}{\sqrt{\pi}}\right),
\quad \Go_1^-(s)=-\fr{2\Gd}{a\Gt_1 s\sqrt{\pi}}.
\label{4.8}
\eeq
Expand the coefficients $A^+_n$ and $B_n^-$  as
\beq
A_n^+=\Gl^{2n+1}\sum_{k=0}^\infty a_{n,k}\Gl^k, \quad B_n^-=\Gl^{2n}\sum_{k=0}^\infty b_{n,k}\Gl^k, 
\label{4.9}
\eeq
and substitute them into the system (\ref{4.7}). This yields
$$
\sum_{k=0}^\infty b_{n,k}\Gl^k=\fr{2}{\pi}\sum_{m=0}^\infty\fr{\Gl^{2m+1}}{n+m+1/2}\sum_{k=0}^\infty a_{m,k}\Gl^k,
$$
\beq
\sum_{k=0}^\infty a_{n,k}\Gl^k=\fr{1}{2\pi}\sum_{m=0}^\infty\fr{\Gl^{2m}}{n+m+1/2}\sum_{k=0}^\infty b_{m,k}\Gl^k-\fr{\Gd^*}{\pi(2n+1)},
\quad n=0,1,\ldots,
\label{4.10}
\eeq
where $\Gd^*=2\Gd(a\Gt_1\sqrt{\pi})^{-1}$.
From here, on comparing the coefficients of the same powers of $\Gl$, we deduce
$$
 a_{n,0}=-\fr{\Gd^*}{2\pi(n+1/2)},  \quad b_{n,0}=0,
 \quad
 a_{n,j}=\fr{b_{0,j}}{2\pi(n+1/2)}, \quad  b_{n,j}=\fr{2a_{0,j-1}}{\pi(n+1/2)},
 $$
 $$
a_{n,j+2}=\fr{1}{2\pi}\left(\fr{b_{0,j+2}}{n+1/2}+\fr{b_{1,j}}{n+3/2}\right), \quad 
b_{n,j+2}=\fr{2}{\pi}\left(\fr{a_{0,j+1}}{n+1/2}+\fr{a_{1,j-1}}{n+3/2}\right),  \quad \ldots, 
$$
$$
a_{n,j+2p}=\fr{1}{2\pi}\sum_{k=0}^p\fr{b_{k,j+2p-2k}}{n+k+1/2}, \quad 
b_{n,j+2p}=\fr{2}{\pi}\sum_{k=0}^p\fr{a_{k,j+2p-2k-1}}{n+k+1/2},\quad 
$$
\beq
n=0,1,\ldots,  \quad j=1,2,  \quad p=0,1,\ldots.
\label{4.11}
\eeq

\setcounter{equation}{0}

\section{Factorization of the triangular matrices. The partial indices of factorization}

In Sections 2 and 3, the vector Riemann--Hilbert problems were solved directly by bypassing factorization of the matrix coefficient $G(s)$.
Here, 
we aim to construct  factorization matrices. This will be done by the method applied in \cite{ant2}  based on the solutions to the homogeneous vector Riemann--Hilbert
problem in an extended class. Similarly to \cite{ant5} we shall  show that these matrices constitute the canonical matrices of factorization and determine the partial indices of factorization. 

\subsection{$2\times 2$ triangular matrix}

Since the solution to the  vector Riemann--Hilbert problem (\ref{4.4}), the functions $\GF_1^\pm(s)$ and  $\GF_2^\pm(s)$, have a fractional order at infinity,
\beq
\GF_1^\pm(s)=O(s^{-3/2}), \quad \GF_2^\pm(s)=O(s^{-1/2}),\quad s\to\infty, \quad s\in \CD^\pm,
\label{4.12}
\eeq
first, we transform the original problem  (\ref{4.4})
into a new one whose solution  has integer orders at infinity. With the aid of the function
\beq 
\tan\fr{\pi s}{2}=L^+(s) L^-(s), \quad s\in\CL,
\label{4.13}
\eeq 
where $L^+(s)$ and $L^-(s)$ are given by (\ref{3.1}), we write
\beq
\left(
\begin{array}{c}
\tilde\GF_1^+(s)\\
\tilde\GF_2^+(s)\\
\end{array}\right)= G_0(s)
\left(
\begin{array}{c}
\tilde\GF_1^-(s)\\
\tilde\GF_2^-(s)\\
\end{array}\right)-
\fr{2\Gl^s\hat w^-(s)L^-(s)}{a\Gt_1}
\left(
\begin{array}{c}
0\\
1\\
\end{array}\right).
\label{4.14}
\eeq
Here,
$$
G_0(s)=
\left(
\begin{array}{cc}
\Gl^{-s}\cot\fr{\pi s}{2} & 0\\
2 & -\Gl^{s+1}\tan\fr{\pi s}{2}\\
\end{array}\right),
$$
\beq
 \tilde\GF_1^+(s)=\fr{\GF_1^+(s)}{L^+(s)}, \quad
 \tilde\GF_2^+(s)=L^+(s)\GF_2^+(s),
\quad
 \tilde\GF_1^-(s)=L^-(s)\GF_1^-(s), \quad  \tilde\GF_2^-(s)=\fr{\GF_2^-(s)}{L^-(s)}.
 \label{4.15}
 \eeq 
Due to (\ref{4.12}) and (\ref{4.15}), the new functions vanish at $s=\infty$ and have an integer order at this point,
$\tilde\GF_j^\pm(s)=O(s^{-1})$,  
$s\in\CD^\pm$,  $s\to\infty,\quad j=1,2.$

We wish to
find two matrices, $X^+(s)$ and $X^-(s)$, analytic in the domains $\CD^+$ and $\CD^-$,
respectively, having a finite order at infinity and solving the following matrix equation:
\beq
X^+(s)=G_0(s)X^-(s), \quad s\in \CL.
\label{4.17}
\eeq
Denote
\beq
X^\pm(s)=\left(\begin{array}{cc}
\Gc_{11}^\pm(s)  & \Gc_{12}^\pm(s)\\
 \Gc_{21}^\pm(s) & \Gc^\pm_{22}(s)\\
 \end{array}
 \right).
 \label{4.18}
 \eeq
On substituting these matrices into (\ref{4.17}) we discover
$$
\Gc_{1l}^+(s)=\Gl^{-s}\cot\fr{\pi s}{2}\Gc_{1l}^-(s),
$$
\beq
\Gc_{2l}^+(s)=2\Gc_{1l}^-(s)-\Gl^{s+1}\tan\fr{\pi s}{2}\Gc_{2l}^-(s), \quad s\in\CL, \quad l=1,2.
\label{4.19}
\eeq
Employing the factorization (\ref{3.1}) of the function $L(s)$, after rearrangement, we arrive at
$$
\Gc_{2l}^+(s)-\Psi_l^+(s)=2\Gc_{1l}^-(s)-\Gl^{s+1}\tan\fr{\pi s}{2}\Gc_{2l}^-(s)-\Psi_l^+(s),
$$
\beq
\Gl \Gc_{2l}^-(s)-\GO_l^-(s)=2\Gc_{1l}^+(s)-\Gl^{-s}\cot\fr{\pi s}{2}\Gc_{2l}^+(s)-\GO_l^-(s).
\label{4.20}
\eeq
Here,
\beq
\Psi_l^+(s)=\sum_{m=0}^\infty\fr{A^+_{lm}}{s-2m-1}, \quad 
\GO_l^-(s)=\sum_{m=0}^\infty\fr{B^-_{lm}}{s+2m}.
\label{4.21}
\eeq
To construct a nontrivial solution, we widen the class of solutions. In the case $l=1$, we choose
$$
\Gc_{21}^+(s)=O(1), \quad \Gc_{11}^+(s)=O(s^{-1}), \quad s\in\CD^+, \quad s\to\infty,
$$
\beq
\Gc_{11}^-(s)=O(1), \quad \Gc_{21}^-(s)=O(s^{-1}), \quad s\in\CD^-, \quad s\to\infty,
\label{4.22}
\eeq
while in the case $l=2$,
$$
\Gc_{22}^+(s)=O(s^{-1}), \quad \Gc_{12}^+(s)=O(1), \quad s\in\CD^+, \quad s\to\infty,
$$
\beq
\Gc_{12}^-(s)=O(s^{-1}), \quad \Gc_{22}^-(s)=O(1), \quad s\in\CD^-, \quad s\to\infty.
\label{4.23}
\eeq
For $l=1$, by the continuity principle and the Liouville theorem, the left- and right-hand sides of the first
equation in (\ref{4.20}) analytically continue each other to the whole complex plane and equal a constant, $C_{11}$.
Without loss, $C_{11}=1$. The second equation gives rise to a constant $C_{12}=0$. Similarly, in the case $l=2$, the corresponding
constants $C_{21}$ (the first equation) and  $C_{22}$ (the second equation) 
have the values $C_{21}=0$ and $C_{22}=1$. On following the procedure described in detail in Section 3 we derive the components of the matrices
of factorization,
the functions $\Gc_{ml}^+(s)$ and  $\Gc_{ml}^-(s)$, in the form
$$
\Gc_{2l}^+(s)=\Psi_l^+(s)+\Gd_{l1},\quad \Gc_{2l}^-(s)=\fr{1}{\Gl}[\GO_l^-(s)+\Gd_{l2}],\quad 
$$$$
\Gc_{1l}^+(s)=\fr12\left\{\GO_l^-(s)+\Gd_{l2}+\Gl^{-s}\cot\fr{\pi s}{2}[\Psi_l^+(s)+\Gd_{l1}]\right\},
$$
\beq
\Gc_{1l}^-(s)=\fr12\left\{\GY_l^+(s)+\Gd_{l1}+\Gl^{s}\tan\fr{\pi s}{2}[\GO_l^-(s)+\Gd_{l2}]\right\}.
\label{4.24}
\eeq
The coefficients $A_{ln}^+$ and $B_{nl}^-$ involved in the representations (\ref{4.21}) of the functions $\Psi^+_l(s)$ and 
$\GO^-_l(s)$ solve the following infinite systems of linear algebraic equations:
$$
A_{ln}^+=\fr{\Gl^{2n+1}}{\pi}\left(\sum_{m=0}^\infty\fr{B_{lm}^-}{n+m+1/2}+2\Gd_{l2}\right),
$$
\beq
B_{ln}^-=\fr{\Gl^{2n}}{\pi}\left(\sum_{m=0}^\infty\fr{A_{lm}^+}{n+m+1/2}-2\Gd_{l1}\right),
\quad n=0,1,\ldots, \quad l=1,2,
\label{4.25}
\eeq
where $\Gd_{lk}$ is the Kronecker symbol, $\Gd_{lk}=1$ if $l=k$ and $0$ otherwise. The solution to the systems (\ref{4.25})
can be represented in the form
\beq
A_{ln}^+=\Gl^{2n+1}\sum_{k=0}^\infty a_{l,n,k}\Gl^k, \quad B_{ln}^-=\Gl^{2n}\sum_{k=0}^\infty b_{l,n,k}\Gl^k
\label{4.26}
\eeq
with the coefficients $a_{l,n,k}$ and $b_{l,n,k}$ being 
recovered from the recurrence relations 
$$
 a_{l,n,0}=\fr{2\Gd_{l2}}{\pi},  \quad b_{l,n,0}=-\fr{2\Gd_{l1}}{\pi},
 \quad
 a_{l,n,j}=\fr{b_{l,0,j}}{\pi(n+1/2)}, \quad  b_{l,n,j}=\fr{a_{l,0,j-1}}{\pi(n+1/2)},
  \quad \ldots, 
$$
$$
a_{l,n,j+2p}=\fr{1}{\pi}\sum_{k=0}^p\fr{b_{l,k,j+2p-2k}}{n+k+1/2}, \quad 
b_{l,n,j+2p}=\fr{1}{\pi}\sum_{k=0}^p\fr{a_{l,k,j+2p-2k-1}}{n+k+1/2},\quad 
$$
\beq
l=1,2, \quad n=0,1,\ldots, \quad j=1,2,
\quad p=0,1,\ldots.
\label{4.27}
\eeq

We have shown that the matrices $X^+(s)$ and $X^-(s)$ with the components (\ref{4.24})  factorize the matrix $G_0(s)$,
$G_0(s)=X^+(s)[X^-(s)]^{-1}$, $s\in\CL$.
We wish to prove next that these matrices constitute 
 the piecewise analytic canonical factorization. 
 We remind that a matrix of factorization is the canonical one if \cite{vek}
 
 (1) $\det X^\pm(s)\ne 0$, $s\in\CD^\pm$, and
 
 (2)   the matrices $X^\pm(s)$ have the normal form at infinity.
 
 A matrix is said to have the normal form  at a point if the order of the determinant at this point is equal to the sum of the orders of the columns.
 The order $\Ga_j$  at $s=\infty$  of a function $y_j(s)$ is determined by   $y_j(s)=\tilde y_j(s)s^{-\Ga_j}$, $s\to\infty$, where the function $\tilde y_j(s)$ is bounded at  infinity
 and $\tilde y_j(\infty)\ne 0$. 
The order $\Ga$ of the vector $\By(s)=(y_1(s),\ldots y_n(s))^T$ at the infinite point
is defined by $\Ga=\min\{\Ga_1,\ldots,\Ga_n\}$. 

Show first that the matrices $X^\pm(s)$ are not singular in any finite part of $\CD^\pm$, that is $\det X^\pm(s)=\Gc_{11}^\pm(s)\Gc_{22}^\pm(s)-\Gc_{12}^\pm(s)\Gc_{21}^\pm(s)\ne 0$. In view of (\ref{4.24}), we have
\beq
\det X^+(s)=-\fr{1}{2}\Gc(s), \quad s\in\CD^+,\quad \det X^-(s)=\fr{\Gc(s)}{2\Gl}, \quad s\in\CD^-,
\label{4.28}
\eeq
where
\beq
\Gc(s)=[1+\Psi^+_1(s)][1+\GO_2^-(s)]-\Psi_2^+(s)\GO_1^-(s).
\label{4.29}
\eeq
The relations (\ref{4.28}) imply that the function $\Gc(s)$ is analytic everywhere in the whole complex plane and $\Gc(s)\sim 1$, $s\to\infty$.
Therefore, $\Gc(s)\equiv 1$ in the whole plane, 
\beq
\det X^+(s)=-\fr{1}{2},\quad s\in\CD^+,\quad \det X^-(s)=\fr{1}{2\Gl}, \quad s\in\CD^-, 
\label{4.30}
\eeq 
 $\det X^\pm(s)\ne 0$ in $\CD^\pm$, and the order of the functions $\det X^\pm(s)$ at $s=\infty$ equals 0.

Analyze next the behavior of the columns of the factorization matrices $X^\pm(s)$, $X^\pm_l(s)=(\Gc_{1l}^\pm(s),\Gc_{2l}^\pm(s))^T$, at infinity.
We have
$$
X_1^+(s)=\left(\begin{array}{c}
s^{-1}\tilde\Gc_{11}^+(s)\\
\Gc_{21}^+(s)\\
\end{array}\right), \quad
X_2^+(s)=\left(\begin{array}{c}
\Gc_{12}^+(s)\\
s^{-1}\tilde\Gc_{22}^+(s)\\
\end{array}\right), 
\quad s\in\CD^+, \quad s\to\infty, 
$$
\beq
X_1^-(s)=\left(\begin{array}{c}
\Gc_{11}^-(s)\\
s^{-1}\tilde\Gc_{21}^-(s)\\
\end{array}\right), \quad
X_2^-(s)=\left(\begin{array}{c}
s^{-1}\tilde\Gc_{12}^-(s)\\
\Gc_{22}^-(s)\\
\end{array}\right), 
\quad s\in\CD^-, \quad s\to\infty.
\label{4.31}
\eeq
Here, $\Gc_{12}^+$,  $\Gc_{21}^+$,  $\Gc_{11}^-$,  $\Gc_{22}^-$,    $\tilde\Gc_{11}^+$,  $\tilde\Gc_{22}^+$, $\tilde\Gc_{12}^-$, and  $\tilde\Gc_{21}^-$ are bounded and nonzero at $s=\infty$.
This implies that the orders at infinity of both of the columns of the matrices $X^+(s)$ and $X^-(s)$ are equal to zero.
According to the definition of the normal form, the matrices $X^\pm(s)$ are normal at infinity.
Since we have also proved that $X^\pm(s)$ are not singular in $\CD^\pm$, we may conclude
that the matrix $X(s)=X^\pm(s)$, $s\in\CD^\pm$, is the canonical matrix of factorization.
The orders of its columns, $\Gk_1=0$ and $\Gk_2=0$, are the partial indices of factorization. 
According to the stability criterion \cite{goh}  applied to an order-2 vector Riemann--Hilbert problem, 
if  $\Gk_1\le \Gk_2$ and $\Gk_2-\Gk_1\le 1$, then the system of partial indices is stable. Thus we conclude that 
the system of partial indices associated with the Riemann--Hilbert problem (\ref{4.14}) is stable.

\subsection{$3\times 3$ triangular matrix}

To deal with functions having the same order-1 at infinity, we employ the diagonal matrix 
\beq
\diag\left\{\cot\fr{\pi s}{2}, \tan\fr{\pi s}{2},  L(s)\right\}=\diag\left\{\fr{1}{L^+(s)L^-(s)},L^+(s)L^-(s),\fr{L^+(s)}{2L^-(s)}\right\} 
\label{4.32}
\eeq
and introduce the new functions
$$
 \tilde\GF_1^+(s)=\fr{\GF_1^+(s)}{L^+(s)}, \quad
 \tilde\GF_2^+(s)=L^+(s)\GF_2^+(s), \quad
 \tilde\GF_3^+(s)=L^+(s)\GF_3^+(s),\quad s\in\CD^+
$$
\beq
 \tilde\GF_1^-(s)=L^-(s)\GF_1^-(s), \quad  \tilde\GF_2^-(s)=\fr{\GF_2^-(s)}{L^-(s)}, \quad  \tilde\GF_3^-(s)=L^-(s)\GF_3^-(s), \quad s\in\CD^-.
\label{4.33}
\eeq
These functions decay at infinity,  $\tilde\GF_j^\pm(s)=O(s^{-1})$, $s\to\infty$, and 
 solve the following Riemann--Hilbert problem:
\beq
\left(
\begin{array}{c}
\tilde\GF_1^+(s)\\
\tilde\GF_2^+(s)\\
\tilde\GF_3^+(s)\\
\end{array}\right)= G_0(s)
\left(
\begin{array}{c}
\tilde\GF_1^-(s)\\
\tilde\GF_2^-(s)\\
\tilde\GF_3^-(s)\\
\end{array}\right)-
\fr{2\Gl_1^s\hat w^-(s)L^-(s)}{a\Gt_1}
\left(
\begin{array}{c}
0\\
0\\
1\\
\end{array}\right).
\label{4.34}
\eeq
where
\beq
G_0(s)=
\left(
\begin{array}{ccc}
\Gl^{-s}\cot\fr{\pi s}{2} & 0  & 0\\
0 & \left(\fr{\Gl_0}{\Gl_1}\right)^{-s-1}\tan\fr{\pi s}{2} & 0 \\ 
2 & -\Gl^{s+1}\tan\fr{\pi s}{2}   & 2\Gl_0^s\\
\end{array}\right).
\label{4.35}
\eeq
With this definition we state the factorization problem $G_0(s)=X^+(s)[X^-(s)]^{-1}$, $s\in\CL$. 
For the entries $\Gc_{ml}^\pm(s)$ $(m,l=1,2,3)$ of the matrices $X^\pm(s)$ we have the system of equations
$$
\Gc_{1l}^+(s)=\Gl_1^{-s}\cot\fr{\pi s}{2}\Gc_{1l}^-(s),
$$
$$
\Gc_{2l}^+(s)=\left(\fr{\Gl_0}{\Gl_1}\right)^{-s-1}\tan\fr{\pi s}{2}\Gc_{2l}^-(s),
$$
\beq
\Gc_{3l}^+(s)=2\Gc_{1l}^-(s)-\Gl_1^{s+1}\tan\fr{\pi s}{2}\Gc_{2l}^-(s)+2\Gl_0^s\Gc_{3l}^-(s), \quad s\in\CL, \quad l=1,2,3.
\label{4.36}
\eeq
Similarly to the $2\times 2$-case, by rearranging the equations, removing the inadmissible poles and extending the class of solutions
by admitting that the properly chosen functions are bounded and nonzero at infinity we arrive at 
$$
\Gc_{3l}^+(s)-\Psi_l^+(s)=2\Gc_{1l}^-(s)-\Gl_1^{s+1}\tan\fr{\pi s}{2}\Gc_{2l}^-(s)+2\Gl_0^s\Gc_{3l}^-(s)-\Psi_l^+(s)=\Gd_{l1},
$$
$$
2\Gc_{1l}^+(s)-\Gl_1^{-s}\cot\fr{\pi s}{2}\Gc_{3l}^+(s)
-\GO_l^+(s)-\GO_l^-(s)
$$$$
=\Gl_1\Gc_{2l}^-(s)-2\left(\fr{\Gl_0}{\Gl_1}\right)^{s}\cot\fr{\pi s}{2}\Gc_{3l}^-(s)-\GO_l^+(s)-\GO_l^-(s)=\Gd_{l2},
$$
\beq
\Gl_0 \Gc_{2l}^+(s)-2\left(\fr{\Gl_0}{\Gl_1}\right)^{-s}\tan\fr{\pi s}{2}\Gc_{1l}^+(s)+\Gl_0^{-s}\Gc_{3l}^+(s)-\Psi_l^-(s)
=2\Gc_{3l}^-(s)-\GY_l^-(s)=\Gd_{l3}.
\label{4.37}
\eeq
Here, $l=1,2,3$, 
$$
\Psi_l^+(s)=\sum_{m=0}^\infty\fr{A^+_{lm}}{s-2m-1}, \quad \Psi_l^-(s)=\sum_{m=0}^\infty\fr{A^-_{lm}}{s+2m+1},
$$
\beq
\GO_l^+(s)=\sum_{m=0}^\infty\fr{B^+_{lm}}{s-2m-2}, \quad \GO_l^-(s)=\sum_{m=0}^\infty\fr{B^-_{lm}}{s+2m}.
\label{4.38}
\eeq
To remove the undesirable poles in (\ref{4.37}), we shall select the coefficients $A_{ln}^\pm$ and $B^\pm_{ln}$ 
as the solution to the following systems of linear algebraic equations: 
$$
B_{ln}^-=\fr{2\Gl_1^{2n}}{\pi}\left(\sum_{m=0}^\infty\fr{A_{lm}^+}{2n+2m+1}-\Gd_{l1}\right),
$$$$
A_{ln}^+=\fr{2\Gl_1^{2n+1}}{\pi}
\left(\sum_{m=0}^\infty\fr{B_{lm}^-}{2n+2m+1}+\sum_{m=0}^\infty\fr{B_{lm}^+}{2n-2m-1}+\Gd_{l2}\right),
$$ 
$$
A_{ln}^-=-\fr{2}{\pi}\left(\fr{\Gl_0}{\Gl_1}\right)^{2n+1}
\left(\sum_{m=0}^\infty\fr{B_{lm}^+}{2n+2m+3}+\sum_{m=0}^\infty\fr{B_{lm}^-}{2n-2m+1}-\Gd_{l2}\right),
$$ 
\beq
B_{ln}^+=-\fr{2}{\pi}\left(\fr{\Gl_0}{\Gl_1}\right)^{2n+2}\left(
\sum_{m=0}^\infty\fr{A_{lm}^-}{2n+2m+3}+\Gd_{l3}\right), \quad n=0,1,\ldots,\quad l=1,2.
\label{4.39}
\eeq
From (\ref{4.37}) we infer that the components of the factorizing matrices have the form
$$
\Gc_{1l}^+(s)=\fr12\left\{\Gl_1^{-s}\cot\fr{\pi s}{2}[\GY_l^+(s)+\Gd_{l1}]+\GO_l^-(s)+\GO_l^+(s)+\Gd_{l2}\right\},
$$
$$
\Gc_{1l}^-(s)=\fr12\left\{\Gl_1^{s}\tan\fr{\pi s}{2}[\GO_l^-(s)+\GO_l^+(s)+\Gd_{l2}]+
\GY_l^+(s)+\Gd_{l1}\right\},
$$
$$
\Gc_{2l}^+(s)=\fr{1}{\Gl_0}\left\{\left(\fr{\Gl_0}{\Gl_1}\right)^{-s}\tan\fr{\pi s}{2}[\GO_l^-(s)+\GO_l^+(s)+\Gd_{l2}]+
\GY_l^-(s)+\Gd_{l3}\right\},
$$
$$
\Gc_{2l}^-(s)=\fr{1}{\Gl_1}\left\{\left(\fr{\Gl_0}{\Gl_1}\right)^{s}\cot\fr{\pi s}{2}[\GY_l^-(s)+\Gd_{l3}]+\GO_l^-(s)+\GO_l^+(s)+\Gd_{l2}
\right\},
$$
\beq
\Gc_{3l}^+(s)=\GY_l^+(s)+\Gd_{l1}, \quad \Gc_{3l}^-(s)=\fr12\left[\GY_l^-(s)+\Gd_{l3}\right].
\label{4.40}
\eeq
Show now that the matrices $X^+(s)$ and $X^-(s)$ are not singular in the domains $\CD^+$ and $\CD^-$, respectively.
By direct computation we obtain
\beq
\det X^+(s)=\fr{\Gc(s)}{2\Gl_0}, \quad  \det X^-(s)=\fr{\Gc(s)}{4\Gl_1},
\label{4.41}
\eeq
where
$$
\Gc(s)=1-[\GO_1^+(s)+\GO_1^-(s)][\GY_2^+(s)+\GY_2^+(s)\GY_3^-(s)-\GY_2^-(s)\GY_3^+(s)]
$$
$$
+[\GO_2^+(s)+\GO_2^-(s)+1][\GY_1^+(s)+\GY_3^-(s)+\GY_1^+(s)\GY_3^-(s)-\GY_1^-(s)\GY_3^+(s)]$$
\beq
-[\GO_3^+(s)+\GO_3^-(s)]
[\GY_2^-(s)+\GY_1^+(s)\GY_2^-(s)-\GY_1^-(s)\GY_2^+(s)]+\GO_2^+(s)+\GO_2^-(s).
\label{4.42}
\eeq
Recall that in Section 4.1, the function $\Gc(s)$ was equal to 1. The same reasoning holds for the function $\Gc(s)$ in the   $3\times 3$ case, 
$\Gc(s)\equiv 1$ in  the whole plane, and
\beq
\det X^+(s)=\fr{1}{2\Gl_0}, \quad  \det X^-(s)=\fr{1}{4\Gl_1}.
\label{4.43}
\eeq
We can immediately conclude that not only the matrices $X^+(s)$ and $X^-(s)$ are not singular in any finite part of the complex $s$-plane but also that
they have zero orders at infinity.
Analyze now the orders of the columns of the matrices $X^+(s)$ and $X^-(s)$. In view of (\ref{4.40}) it is seen that
two elements of each columns have order 1, while the third entry has order 0. Therefore the orders of all columns equal 0.
The matrices $X^+(s)$ and $X^-(s)$ are normal at infinity, not singular everywhere in the domains $\CD^+$ and $\CD^-$ and therefore
they constitute the canonical factorization of the matrix coefficient $G_0(s)$ of the Riemann--Hilbert problem (\ref{4.34}).
Since the orders at the infinite point of the columns of these matrices are zeros, the partial indices of factorization, $\Gk_1$, $\Gk_2$, and $\Gk_3$, are also equal to zero.  

\setcounter{equation}{0}

\section{Contact stresses and normal displacements in the case of a circular inclusion. Numerical results}

Suppose that a rigid inclusion and a crack are both penny-shaped, the inclusion is flat, $w(r)=\Gd>0$,  $0\le r\le b$, and it is
is planted between the crack faces. In this case the contact area is known, $0\le r\le b$, $\Gl=b/a<1$, and the Mellin transforms (\ref{4.3}) of the contact stresses
$\Gs_z=\psi_0(r)$, $0\le r\le b$,  and $\Gs_z=\psi_1(r)$, $r\ge a$, and the normal displacement $\Gc_1(r)$ in the annulus $b\le r\le a$
have been found. They are given by (\ref{4.6}).
To derive the contact stresses (the normal traction) and the normal displacements we need to invert the Melin transforms and
rewrite the integrals in the form convenient for computations.  We have 
\beq
\Gs_z(r,0^+)=\fr{1}{2\pi i\Gl}\int_\CL L^-(s)\left[\GO^-(s)-\fr{\Gd^*}{s}\right]\left(\fr{r}{b}\right)^{-s-1}ds, \quad 0\le r<b.
\label{5.1}
\eeq
where $\Gd^*=2\Gd(a\Gt_1\sqrt{\pi})^{-1}$. On employing the residues theory and changing the order of summation we may eventually write for $0<r<b$
\beq
\Gs_z(r,0^+)=-\fr{\Gd^*}{\Gl\sqrt{\pi[1-(r/b)^2]}}-\fr{1}{2\Gl\sqrt{\pi}}
\sum_{m=0}^\infty\fr{B_m^-}{m-1/2}F\left(\fr32,\fr12-m; \fr32-m;\fr{r^2}{b^2}\right).
\label{5.2}
\eeq
The series converges rapidly due to the exponential decay of the coefficients $B_m^-$ as $m\to\infty$. Now, if $r\to b$, it is convenient
to use formula 9.131(2) \cite{gra} to obtain
\beq
\Gs_z(r,0^+)=\fr{1}{\Gl\sqrt{\pi[1-(r/b)^2]}}\left[-\Gd^*+
\sum_{m=0}^\infty B_m^-\sum_{j=0}^m
\fr{(-m)_j}{(1/2)_j}\left(1-\fr{r^2}{b^2}\right)^j\right],
\label{5.3}
\eeq
where $0<r<b$ and $(a)_m=a(a+1)\ldots(a+m-1)$ is the factorial symbol.
\begin{figure}[t]
\centerline{
\scalebox{0.5}{\includegraphics{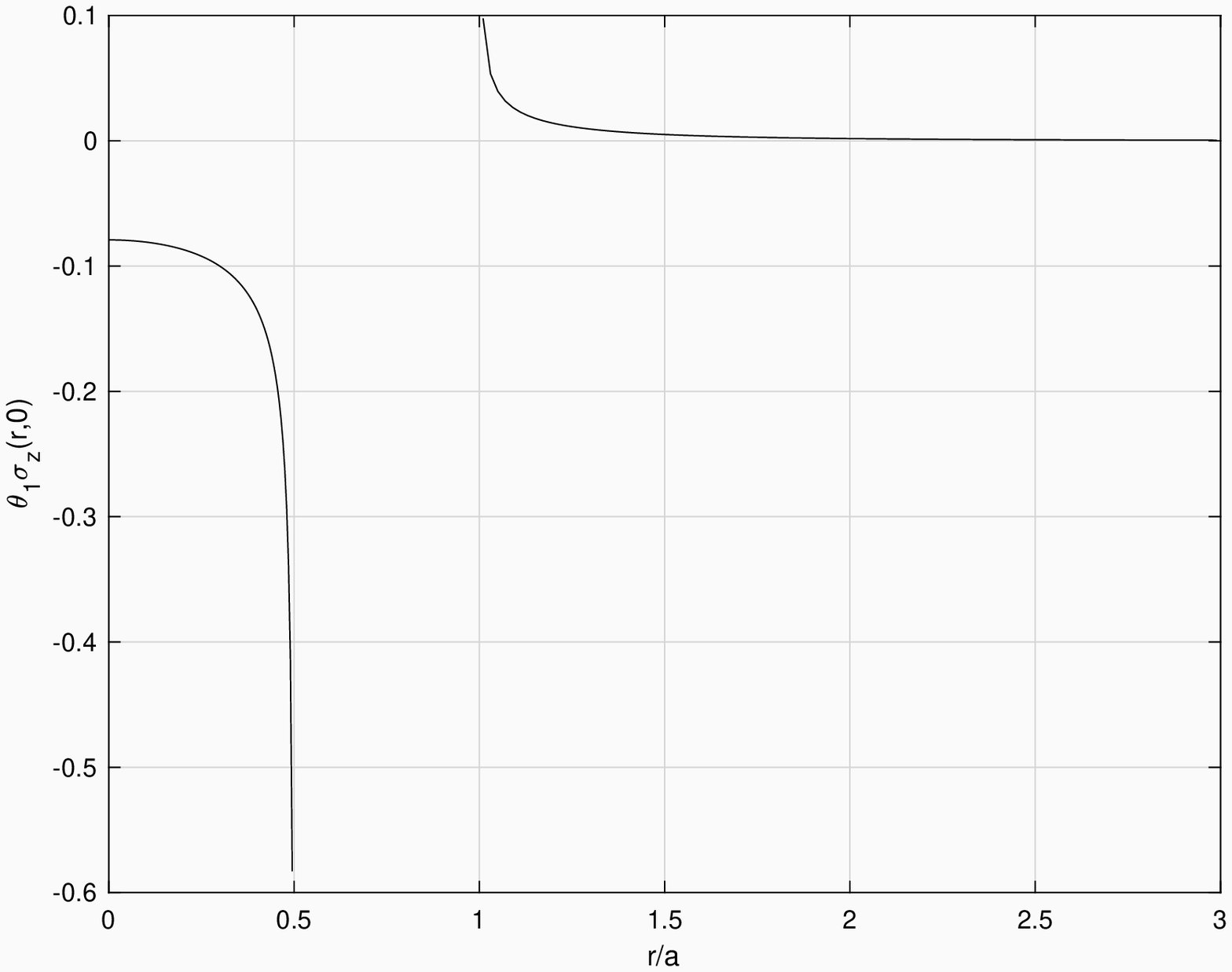}}}
\caption{Normal stress $\Gt_1\Gs_z(r,0)$ for $0\le r<b$ and $r>a$ when $\Gl=0.5$ and $\Gd/a=0.05$.} 
\label{fig1}
\end{figure}

Let now $a<r<\infty$. By inversion of the Mellin integral $\GF_2^+(s)$ and using its representation (\ref{4.6}) we find
\beq
\Gs_z(r,0)=\fr{1}{\pi i}\int_\CL\fr{\GY^+(s)}{L^+(s)}\left(\fr{r}{a}\right)^{-s-1}ds, \quad r>a.
\label{5.4}
\eeq
Similarly to the integral (\ref{5.1}) we deduce the series representation in terms of the Gauss function
\beq
\Gs_z(r,0)=\fr{1}{\sqrt{\pi}}\left(\fr{a}{r}\right)^3
\sum_{m=0}^\infty\fr{A_m^+}{m-1/2}F\left(\fr32,\fr12-m; \fr32-m;\fr{a^2}{r^2}\right), \quad r>a.
\label{5.5}
\eeq
In the contact zone $r>a$ when $r$ is close to $a$ this formula can be rewritten in the form
\beq
\Gs_z(r,0)=-\fr{2}{\sqrt{\pi[1-(a/r)^2]}}\left(\fr{a}{r}\right)^3
\sum_{m=0}^\infty A_m^+
\sum_{j=0}^m
\fr{(-m)_j}{(1/2)_j}\left(1-\fr{a^2}{r^2}\right)^j,\quad r>a.
\label{5.6}
\eeq

\begin{figure}[t]
\centerline{
\scalebox{0.5}{\includegraphics{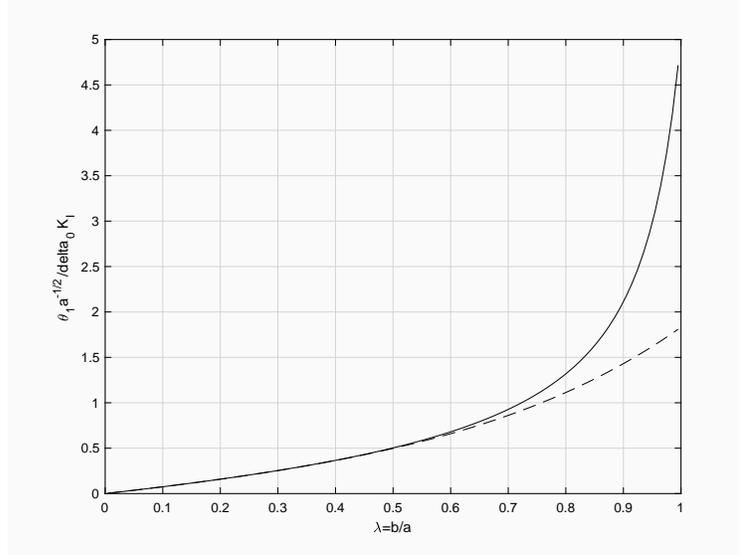}}}
\caption{Normlized stress intensity factor $\Gt_1a^{-1/2}K_I/\Gd_0$ at the crack  tip $r=a$  versus the parameter $\Gl$: the solid line corresponds to the exact formula (\ref{5.8}), while
the dashed line relates to the asymptotic formula (\ref{5.8.3}).} 
\label{fig2}
\end{figure} 
Formulas (\ref{5.3}) and (\ref{5.6}) indicate that the normal stresses have the square root singularity as $r\to b^-$ and $r\to a^+$. This is consistent with
the graphs of $\Gs_z(r,0)$ in the contact zone $0\le r<b$ and for $r>a$ as $z=0$ (Fig.1). From the last formula 
we may immediately derive the stress intensity factor at the tip $r=a$ of the crack
\beq
K_I=\lim_{r\to a^+}\sqrt{2\pi(r-a)}\Gs_z(r,0).
\label{5.7}
\eeq
It is given by
\beq
K_I=-2\sqrt{a}\sum_{m=0}^\infty A_m^+.
\label{5.8}
\eeq
Note that the same formula is obtained directly from the expression (\ref{4.6}) for the integral $\GF_2^+(s)$ my
employing the Abelian theorems for the Mellin transforms. In addition to the exact formula (\ref{5.8}), it is possible
to write a simple asymptotic formula in terms of $\Gl^j$. By virtue of the first formula in (\ref{4.9}) and (\ref{5.8}) we have
\beq
K_I=-2\sqrt{a}\Gl[a_{00}+\Gl a_{01}+\Gl^2(a_{02}+a_{10})+\Gl^3(a_{03}+a_{11})+\Gl^4(a_{04}+a_{12}+a_{20})+\ldots].
\label{5.8.1}
\eeq
We next employ the recurrence relations (\ref{4.11}) and deduce the expressions
$$
a_{00}=-\fr{\Gd^*}{\pi}, \quad a_{10}=-\fr{\Gd^*}{3\pi}, \quad  a_{20}=-\fr{\Gd^*}{5\pi},
$$$$
 a_{01}=-\fr{4\Gd^*}{\pi^3},\quad  a_{02}=-\fr{16\Gd^*}{\pi^5},\quad  a_{11}=-\fr{4\Gd^*}{3\pi^3},\quad  a_{12}=-\fr{16\Gd^*}{3\pi^5},
$$
\beq
 a_{03}=-\fr{4\Gd^*}{\pi^3}\left(\fr{16}{\pi^4}+\fr29\right),\quad  a_{04}=-\fr{64\Gd^*}{\pi^5}\left(\fr{4}{\pi^4}+\fr19\right).
\label{5.8.2}
\eeq
Here, as before, $\Gd^*=2\Gd(a\Gt_1\sqrt{\pi})^{-1}$. Substituting these formulas into (\ref{5.8.1}) yields the asymptotic expansion of the
coefficient $K_I$ for small $\Gl$ 
$$
K_I=\fr{4\sqrt{a}\Gd_0}{\pi^{3/2}\Gt_1}
\left[\Gl+\fr{4\Gl^2}{\pi^2}+\left(\fr{16}{\pi^4}+\fr13\right)\Gl^3
\right.
$$
\beq
\left.
+\fr{4}{\pi^2}\left(\fr{16}{\pi^4}+\fr59\right)\Gl^4+\left(\fr{256}{\pi^8}+\fr{112}{9\pi^4}+\fr15\right)\Gl^5+O(\Gl^6)\right],
\label{5.8.3}
\eeq
where $\Gd_0=\Gd/a$. Referring to Fig. 2 we conclude that for $0<\Gl<0.6$ the asymptotic expansion (\ref{5.8.3}) is in good agreement with the exact formula
(\ref{5.8}). 

\begin{figure}[t]
\centerline{
\scalebox{0.5}{\includegraphics{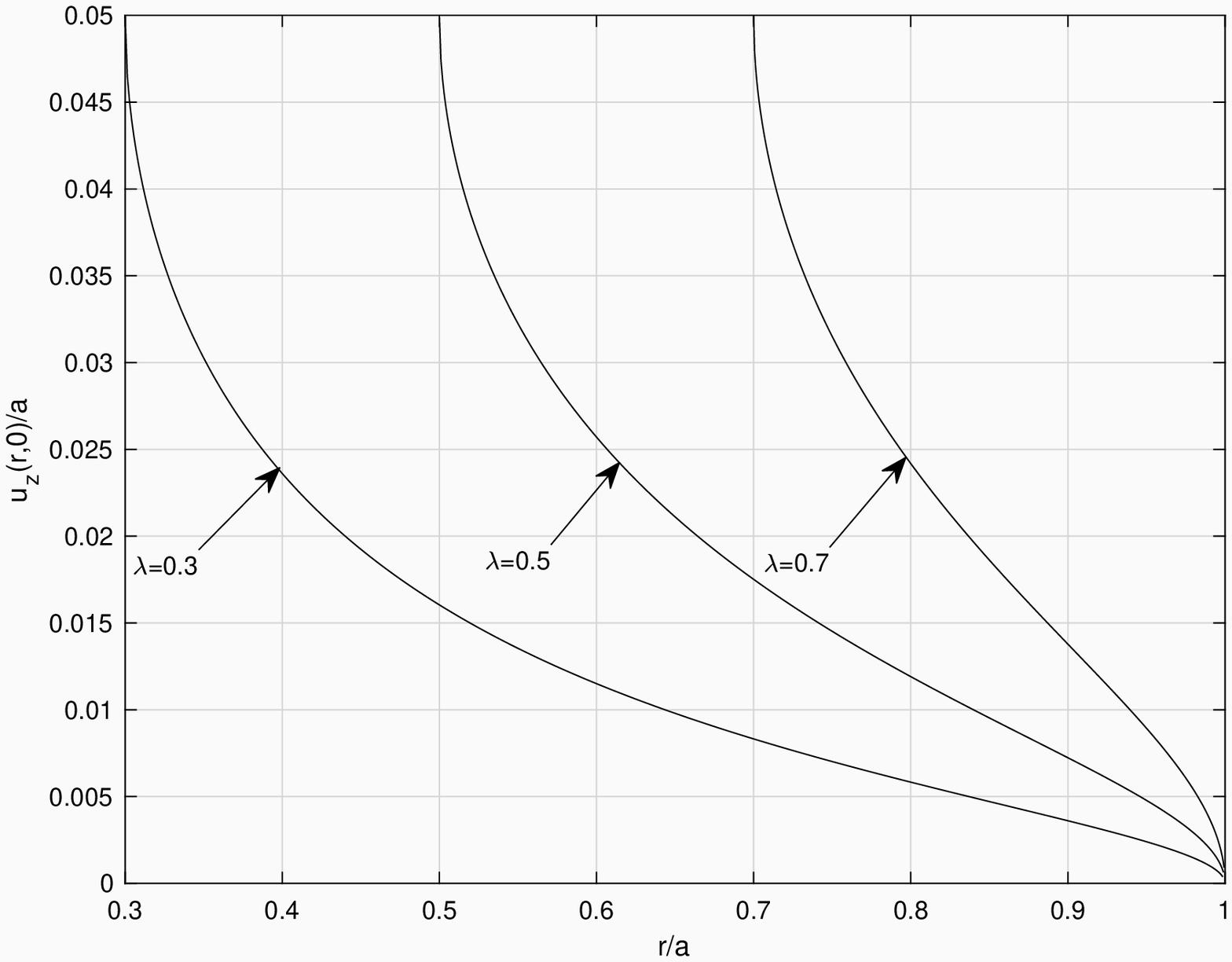}}}
\caption{Normal displacement $u_z(r,0)/a$ for $b<r<a$ when $\Gl=0.3$, $\Gl=0.5$, and $\Gl=0.7$ for $\Gd/a=0.05$.} 
\label{fig3}
\end{figure} 
The profile of the crack annular surface is described by the function $u_z(r,0^+)=-a\Gt_1\Gc_1(r)$, $b\le r\le a$.  
Its integral representation is derived by the Mellin inversion of the function $\GF_1^-(s)$ given by (\ref{4.6}) that is
$$
\Gc_1(r)=\fr{1}{2\pi i}\int_\CL\fr{\GY^+(s)}{L^-(s)}\left(\fr{r}{a}\right)^{-s}ds
$$
\beq
+\fr{1}{4\pi i}\int_\CL\left\{L^+(s)\left[\GO^-(s)-\fr{\Gd^*}{s}\right]+\fr{\sqrt{\pi}\Gd^*}{s}\right\}\left(\fr{r}{b}\right)^{-s}ds,
\quad b<r<a.
\label{5.9}
\eeq
As before we employ the theory of residues and, in addition, formula 9.121(26) from \cite{gra}
\beq
F\left(\fr12,\fr12;\fr32;x^2\right)=\fr{\sin^{-1} x}{x}.
\label{5.9'}
\eeq
 This gives rise to the 
series representation for $b<r<a$
\beq
\Gc_1(r)=-\fr{2\Gd}{\pi a\Gt_1}\sin^{-1}\fr{b}{r}+\fr{b}{\sqrt{\pi}r}\sum_{m=0}^\infty \fr{B^-_m}{2m+1}f_m\left(\fr{b^2}{r^2}\right)-
\fr{2}{\sqrt{\pi}}\sum_{m=0}^\infty \fr{A^+_m}{2m+1}f_m\left(\fr{r^2}{a^2}\right),
\label{5.10}
\eeq
where
\beq
f_m(x)=F\left(\fr12,m+\fr12; m+\fr32;x\right).
\label{5.11}
\eeq
For according to formula 9.131(2)  from \cite{gra}  we can represent the function $f_m(x)$ for $x$ close to $1^-$ as
\beq
f_m(x)=\fr{\sqrt{\pi}(m+1/2)}{m!}
\sum_{j=0}^\infty
\left[\fr{\GG(m+j+1/2)}{j!}-\fr{\GG(m+j+1)\sqrt{1-x}}{\GG(j+3/2)}\right](1-x)^j,
\label{5.12}
\eeq
 and, in the limit,
 \beq
 f_m(1^-)=\fr{\pi(3/2)_m}{2m!}.
 \label{5.13}
 \eeq
We wish now to verify that the normal displacement is continuous at the points $r=b$ and $r=a$.
In view of (\ref{5.10}) and (\ref{5.13}) we deduce
\beq
\Gc_1(b^+)=-\fr{\Gd}{a\Gt_1}+\fr12
\sum_{m=0}^\infty\fr{B_m^-\GG(m+1/2)}{m!}-
\fr{2}{\sqrt{\pi}}\sum_{m=0}^\infty
\fr{A_m^+ f_m(\Gl^2)}{2m+1}
\label{5.14}
\eeq
and
\beq
\Gc_1(a^-)=-\fr{2\Gd}{\pi a\Gt_1}\sin^{-1}\Gl +\fr{\Gl}{\sqrt{\pi}}
\sum_{m=0}^\infty\fr{B_m^-f_m(\Gl^2)}{2m+1}-
\sum_{m=0}^\infty
\fr{A_m^+\GG(m+1/2)}{m!}.
\label{5.15}
\eeq
Analyze the expression (\ref{5.14}) first. 
On employing formula (\ref{5.11}), changing the order of summation in the second term in (\ref{5.14}) we arrive at
\beq
\Gc_1(b^+)=-\fr{\Gd}{a\Gt_1}+\fr12
\sum_{m=0}^\infty\fr{\GG(m+1/2)}{m!}\left(
B_m^-
-\fr{2\Gl^{2m}}{\pi}\sum_{j=0}^\infty\fr{A_j^+}{m+j+1/2}
\right).
\label{5.16}
\eeq
Due to the first equation in (\ref{4.7}) the expression in the brackets is equal to zero and therefore $u_z(b^+)=-a\Gt_1\Gc_1(b^+)=\Gd$,
and the normal displacement is continuous at $r=b$.

Similarly, formula (\ref{5.15}) becomes
\beq
\Gc_1(a^-)=-\fr{2\Gd\sin^{-1}\Gl}{\pi a\Gt_1}-\sum_{m=0}^\infty\fr{\GG(m+1/2)}{m!}\left(A_m^+-\fr{\Gl^{2m+1}}{2\pi}\sum_{j=0}^\infty\fr{B_j^-}{m+j+1/2}\right). 
\label{5.19}
\eeq
Now it is turn of the second equation in the system (\ref{4.7}). If, in addition, formulas (\ref{4.8}) and  (\ref{5.9'}) are used, then we have $\Gc_1(a^-)=0$, and
the normal displacement is continuous at the point $r=a$ as well. The normalized displacements $u_z(r,0)/a$ for $\Gl=0.3$, $\Gl=0.5$, and $\Gl=0.7$ are plotted
in Fig. 3.

\section{Conclusion}

We  developed  an analytical solution to two model problems of a penny-shaped crack when an annulus-shaped (model 1) or a disc-shaped (model 2) rigid inclusion planted between the crack faces.
The method we proposed for these models recast the governing integral equations with the Weber--Sonin kernel
on two segments as vector Riemann--Hilbert problems with a $3\times 3$ and
$2\times 2$ triangular matrix coefficient.  The solution presented
for model 2 
may be classified as an exact solution since it is given in terms of explicitly defined functions
and exponentially convergent series whose coefficients are defined explicitly in terms of certain recurrence relations.  Similar relations can be also written for model 1.
For model  2, we derived representation formulas for the normal stress and displacement.
For the stress intensity factor, in addition to the exact formula, we gave
a simple asymptotic expansion in terms of $(b/a)^n$, $a$ and $b$ are the crack and inclusion radii, respectively, and $b<a$.
For both models, we also found the canonical matrix of factorization and the partial indices of
factorization which turn out to be zeros and therefore stable.

 \vspace{.1in}

{\bf Data accessibility.}
 No software generated data were created during this study.
 
{\bf Competing interests.} We have no competing interests.

{\bf Funding statement.}
YAA  thanks the Isaac Newton Institute for Mathematical Sciences, Cambridge, for support and hospitality during the programme Complex analysis: techniques, applications and computations,  where a part of work on this paper was undertaken. This work was supported by EPSRC grant no EP/R014604/1 and the Simons Foundation.


\begin{thebibliography}{99}

\bibitem{sne1}  Sneddon IN. 1946
The distribution of stress in the neighbourhood of a crack in an elastic solid.  \textit{Proc. R. Soc.  A }  \textbf{187}, 229-260. 


\bibitem{mos1} Mossakovskii   VI. 1954  A fundamental mixed problem of the theory of elasticity for a half-space with a circular curve of separation of the boundary conditions. 
 \textit{Prikl. Mat. Meh.}  \textbf{18} (1954), 187-196. 


\bibitem{mos2} Mossakovskii  VI, Rybka MT. 1964 Generalization of the Griffith-Sneddon criterion for the case of a nonhomogeneous body.
 \textit{J. Appl. Math. Mech. }  \textbf{28} (1964),  1277-1286.

\bibitem{sne2} 
Sneddon IN,  Lowengrub M. 1969
 \textit{Crack problems in the classical theory of elasticity.} New York: John Wiley \& Sons.

\bibitem{wil} 
Willis JR.  1972 The penny-shaped crack on an interface.  \textit{Quart. J. Mech. Appl. Math.}
 \textbf{25}(3). (1972), 367-385.

\bibitem{ant1}   Antipov YA,  Mkhitaryan SM. 2020
Correspondence principle in plane and axisymmetric mixed boundary-value problems of elasticity.   \textit{Quart. Appl. Math.} Published electronically on June 20 2019.

\bibitem{sel}  Selvadurai APS,  Singh BM. 1984 On the expansion of a penny-shaped crack by a rigid circular disc inclusion.  \textit{Int. J. Fracture}   \textbf{25}, 69-77.

\bibitem{ant2}   Antipov YA. 1987 Exact solution of the problem of pressing an annular stamp into a half-space. 
 \textit{Dokl. Akad. Nauk Ukrain. SSR Ser. A} \textbf{7}, 29-33.

\bibitem{ant3}   Antipov YA. 2015 Vector Riemann-Hilbert problem with almost periodic and meromorphic coefficients and applications.  \textit{Proc. A.} \textbf{471}, no. 2180, 
20150262,

\bibitem{ant}   Antipov YA,  Mkhitaryan SM. 2017 A crack induced by a thin rigid inclusion partly debonded from the matrix,  \textit{Quart. J. Mech. Appl. Math.} \textbf{70},  153-185.

\bibitem{gra} Gradshte\u in IS,  Ryzhik IM. 2007  \textit{Table of Integrals, Series and Products}. 
Oxford: Academic Press.

\bibitem{ant4}   Antipov YA, Popov GYa, Yatsko SI. 1987 Solution of the problem of stress concentration around intersecting defects by using the Riemann problem with an infinite index. 
\textit{J. Appl. Math. Mech.}  \textbf{51}, 357-365.

\bibitem{ant5}   Antipov YA, Silvestrov VV. 2002 Factorization on a Riemann surface in scattering theory. \textit{Quart. J. Mech. Appl. Math.} \textbf{55},  607-654.

\bibitem{vek} Vekua NP.1967  \textit{Systems of Singular Integral Equations}.
 Groningen: Noordhoff.

\bibitem{goh}  Gohberg IC, Krein MG. 1958 On the stability of a system of partial indices of the Hilbert problem for several unknown functions. \textit{Dokl. AN SSSR}. \textbf{119}, 854-857.


\end{thebibliography}
\end{document}